\documentclass{svproc}
\usepackage{url,xcolor}

\usepackage[utf8]{inputenc}
\usepackage{amssymb,amsmath,amsfonts,mathtools,hyperref}
\def\su{{\mathfrak{su}}}
\def\os{{\mathfrak{o}}}

\begin{document}
\mainmatter              
\title{Howe duality and algebras of the Askey--Wilson type: an overview}
\titlerunning{Howe duality and AW algebras}  
%
\author{Julien Gaboriaud \inst{1,2}, Luc Vinet \inst{1, 2} and St\'ephane Vinet \inst{2}}
\authorrunning{J. Gaboriaud et al.} 
\tocauthor{Luc Vinet}
\institute{Centre de Recherches Math\'ematiques, \\
\and D\'epartement de Physique, Universit\'e de Montr\'eal, \\
\small~P.O. Box 6128, Centre-ville Station, Montr\'eal (Qu\'ebec), H3C 3J7, Canada.
\email{julien.gaboriaud@UMontreal.CA}\\
\email{vinet@CRM.UMontreal.CA}\\
\email{stephane.vinet@UMontreal.CA}}
\maketitle

\begin{abstract}
The Askey--Wilson algebra and its relatives such as the Racah and Bannai--Ito algebras were initially
introduced in connection with the eponym orthogonal polynomials. They have since proved ubiquitous. In
particular they admit presentations as commutants that are related through Howe duality. This paper surveys
these results.
\end{abstract}

\section{Introduction}\label{sec:intro}

The quadratic algebras of Askey--Wilson type such as the Askey--Wilson algebra itself, the Racah and
Bannai--Ito algebras and their specializations and contractions encode the bispectral properties of orthogonal
polynomials that arise in recoupling coefficients such as the Clebsch--Gordan or Racah coefficients. It is
therefore natural that these algebras be encountered as centralizers of the diagonal action of an algebra of
interest $\mathfrak{g'}$ such as $\mathfrak{sl}(2)$, $\mathfrak{osp}(1|2)$ or $U_q(\mathfrak{sl}(2))$, on
$n$-fold tensor products of representations of $\mathfrak{g'}$. Indeed, elements of these centralizers will be
used as labelling operators to define bases whose overlaps will be expressed in terms of the corresponding
orthogonal polynomials.

Often the algebra $\mathfrak{g'}$ forms a reductive pair with another algebra $\mathfrak{g}$ in which
case the Howe duality operates in certain modules. This leads to alternative characterizations of the
quadratic algebras that are in correspondance: on the one hand commutants in representations of the universal
enveloping algebra $U(\mathfrak{g})$ and on the other hand, realizations of the type mentioned above as
centralizers in recoupling problems for $\mathfrak{g'}$. This is the topic of this brief review which is
organized as follows. Section \ref{sec:general} presents the general framework. Section \ref{sec:racah}
describes as illustration the dual commutant picture for the Racah algebra; this will involve the reductive
pair $(\mathfrak{o}(6), \mathfrak{su}(1,1))$. Section \ref{sec:more} gives a summary of the different cases
that have been analyzed and Section \ref{sec:conclusion} provides a short outlook.

\section{General Framework}\label{sec:general}

We shall say following \cite{RCR} that two algebras $\mathfrak{g}$ and $\mathfrak{g'}$ have dual
representations on a Hilbert space $\mathcal{H}$ if (1) this space carries fully reducible representations of
both $\mathfrak{g}$ and $\mathfrak{g'}$, (2) the action of $\mathfrak{g}$ and $\mathfrak{g'}$ commute, (3) the
representation $\rho$ of the direct sum $\mathfrak{g} \oplus \mathfrak{g'}$ defined by the actions of
$\mathfrak{g}$ and $\mathfrak{g'}$ on $\mathcal{H}$ is multiplicity-free and (4) each irreducible
representation of $\mathfrak{g}$ occurring in the decomposition of $\rho$ is paired with a unique irreducible
representation of $\mathfrak{g'}$ and vice-versa. This is the essence of Howe duality which can be proved in a
number of situations. We shall consider such instances in this paper.

Consider now a setup with the representation of $\mathfrak{g'}$ in $\mathcal{H} = V^{\otimes2n}$ given by
$\bar{\sigma}^{\otimes{2n}}[\Delta^{(2n-1)}(\mathfrak{g'})]$ where $\bar{\sigma} : \mathfrak{g'} \rightarrow
End \;V$ is a representation of $\mathfrak{g'}$ on the vector space $V$,~ $\Delta : \mathfrak{g'} \rightarrow
\mathfrak{g'} \otimes \mathfrak{g'}$ is the coproduct and $\Delta^{(n)}$ is defined recursively by
$\Delta^{(n)} = (\Delta \otimes 1^{\otimes (n-1)}) \circ \Delta^{(n-1)}$, with $\Delta^{(0)}=1$. This
symmetric situation makes it natural that there be an action of some other algebra $\mathfrak{g}$ on the
carrier space $\mathcal{H}$ that commutes with the action of $\mathfrak{g'}$.  Take the maximal Abelian
subalgebra $\mathfrak{h}$ of $\mathfrak{g}$ to be $\mathfrak{h} \simeq \mathfrak{X}^{\oplus n}$ with
$\mathfrak{X}$ one-dimensional. The pairing under Howe duality with the representations of
$\mathfrak{X}^{\oplus n}$ implies that $\bar{\sigma}^{\otimes{2n}}[\Delta^{(2n-1)}(\mathfrak{g'})] =
\bar{\sigma}^{\otimes{2n}}[\Delta^{\otimes n} \circ \Delta^{(n-1)}(\mathfrak{g'})]$ decomposes into
representations of the form  $\sigma_1 \otimes \sigma_2 \otimes \dots \otimes \sigma_n
(\Delta^{(n-1)}(\mathfrak{g'}))$ with the $\sigma_i$'s being irreducible representations arising in the
decomposition of $\bar{\sigma}^{\otimes2}$. This quotienting by $\mathfrak{h} $ is a way of posing a
generalized Racah problem for the recoupling of the $n$ representations $\sigma_i$ of $\mathfrak{g'}$.

We indicated in the introduction that the quadratic algebras $\mathcal{A}$ of Askey--Wilson type can be
obtained as (subalgebras of) centralizers of diagonal actions in $n$-fold tensor products of representations.
The intermediate Casimir elements in $\sigma_1 \otimes \sigma_2 \otimes \dots \otimes \sigma_n$ manifestly
centralize the action of $\mathfrak{g'}$ on $\mathcal{H}\mod\mathfrak{h}$. They are taken to generate the
quadratic algebra of interest. This provides the first presentation of $\mathcal{A}$ as a commutant. The dual
one is identified as follows in the present context. We know that  $\mathfrak{g}$ is the commutant of
$\mathfrak{g'}$ in $\mathcal{H}$. Moreover from the application of Howe duality, the generators of the
representation $\sigma_1 \otimes \sigma_2 \otimes \dots \otimes \sigma_n$ of $\mathfrak{g'}$ are known to
commute with those that represent the subalgebra $\mathfrak{h} \simeq \mathfrak{X}^{\oplus n}$. The
non-trivial part of the centralizer of $\sigma_1 \otimes \sigma_2 \otimes \dots \otimes \sigma_n$ must
therefore be obtained, in the given representation on $\mathcal{H}\mod\mathfrak{h}$, by those elements in the
universal enveloping algebra of $\mathfrak{g}$ that commute with $\mathfrak{X}^{\oplus n}$. In other words,
$\mathcal{A}$ can also be identified as the commutant of $\mathfrak{h}  \subset \mathfrak{g}$ in
$\mathcal{U}(\mathfrak{g})$ as represented on $\mathcal{H}$.

There is an equivalent way of looking at this. The pairing of the representations of $\mathfrak{g}$ and
$\mathfrak{g'}$ through Howe duality manifests itself in the fact that the Casimir elements of $\mathfrak{g}$
and  $\mathfrak{g'}$ are affinely related. Let $\mathcal{C}$ be a Casimir element of $\mathfrak{g'}$.
Consider for example the intermediate Casimir element given by $\bar{\sigma}^{\otimes{4}}[((\Delta \otimes
\Delta) \circ \Delta) (\mathcal{C})] \otimes 1^{\otimes{(2n-4)}}$ corresponding to the embedding of
$\mathfrak{g'}$ in the first four factors of $\mathfrak{g'}^{ \otimes 2n}$. There will be a subalgebra
$\mathfrak{g}_1$ of $\mathfrak{g}$ that will be dually related to $\mathfrak{g'}$ on the restriction of
$\mathcal{H}$ to $V^{\otimes 4}$ so that its Casimir element will be essentially the one of $\mathfrak{g'}$.
Next, looking at the intermediate Casimir element of $\mathfrak{g'}$ associated to a different embedding, for
instance in the four last factors of $\mathfrak{g'}^{\otimes 2n}$, there will be a dual pairing with a
different embedding in  $\mathfrak{g}$ of the same subalgebra $\mathfrak{g}_1$ and again the two Casimir
elements will basically coincide. These observations lead to the conclusion that the set of intermediate
Casimir elements associated to the representation of $\mathfrak{g'}$ is algebraically identical to the set of
Casimir elements of the subalgebras of $\mathfrak{g}$ that form dual pairs with $\mathfrak{g'}$ when
intermediate representations of the latter are taken. It is not difficult to convince oneself that the set of
invariants connected to the relevant subalgebras of $\mathfrak{g}$ consists in the commutant of the maximal
Abelian subalgebra of $\mathfrak{g}$ as concluded differently before.

To summarize, in situations where Howe duality prevails with $(\mathfrak{g}, \mathfrak{g'})$ the pair of
algebras that are dually represented on $\mathcal{H}$ and if the representation of $\mathfrak{g'}$ is of the
form $\bar{\sigma}^{\otimes{2n}}[\Delta^{(2n-1)}(\mathfrak{g'})]$, the quadratic algebras $\mathcal{A}$ of
Askey--Wilson type can be viewed on one hand as the commutant of this action of $\mathfrak{g'}$ on
$\mathcal{H}$ and thus generated by the intermediate Casimir elements of $\mathfrak{g'}$, or on the other hand
as the commutant of $\mathfrak{h} \subset \mathfrak{g}$ in the intervening representation of
$\mathcal{U}(\mathfrak{g})$. We shall present next an example of how this can be concretely realized.

\section{The dual presentations of the Racah algebra}\label{sec:racah}

The Racah algebra $\mathcal{R}$ has three generators $K_1$, $K_2$, $K_3 $ that are subjected to the relations
\cite{GVZ}:
\begin{align}\label{CR}
\begin{aligned}{}
 [K_1,K_2]=K_3, \qquad\quad\  [K_2,K_3]&={K_2}^2+\{K_1,K_2\}+dK_2+e_1,\\
 [K_3,K_1]&={K_1}^2+\{K_1,K_2\}+dK_1+e_2,
\end{aligned}
\end{align}
where $[A,B] = AB-BA$,~ $\{A,B\}=AB+BA$ and $d$, $e_1$, $e_2$ are central.

We shall explain how dual presentations of the algebra $\mathcal{R}$ as a commutant are obtained in the
fashion described in Section 2. The dual pair will be $(\mathfrak{o}(6), \mathfrak{su}(1,1))$ and the
representation space $\mathcal{H}$ will be that of the state space of six quantum harmonic oscillators with
annihilation and creation operators $a_{\mu} , a_{\nu}^{\dagger},~ \mu, \nu = 1, \dots, 6$ verifying $[a_{\mu}
, a_{\nu}^{\dagger}] = \delta_{\mu\nu}$. The corresponding Hamiltonian ${H = a_1^{\dagger}a_1 + \dots +
a_6^{\dagger}a_6}$ is manifestly invariant under the rotations in six dimensions. These are encoded in the Lie
algebra $\mathfrak{o}(6)$, realized by the generators $L_{\mu\nu} = a_{\mu}^{\dagger}a_{\nu} -
a_{\mu}a_{\nu}^{\dagger}$ and possessing the Casimir element $\mathcal{C}=\sum_{\mu<\nu}{L_{\mu\nu}}^2$.

The Lie algebra $\su(1,1)$ has generators $J_0, J_\pm$ that obey the following commutation relations: $
[J_0,J_\pm]=\pm J_\pm,~ [J_+,J_-]=-2J_0$, and its Casimir operator is given by $C={J_0}^2-J_+J_--J_0$.
The six harmonic oscillators also provide a realization of this algebra through the addition of six
copies of the metaplectic representation of $\mathfrak{su}(1,1)$, for which the generators are mapped to:
$J_0^{(\mu)}=\frac{1}{2}(a_{\mu}^{\dagger}a_{\mu} +\frac{1}{2})$, ${J_+^{(\mu)} =
\frac{1}{2}(a_{\mu}^{\dagger})^2}$, ${J_-^{(\mu)}=\frac{1}{2}(a_{\mu})^2}$, $\mu = 1, \dots, 6$. Note that
the operators $ \sum_{\mu = 1}^{6} J_\bullet^{(\mu)}$ are invariant under rotations.  The space of state
vectors $\mathcal{H}$ thus carries commuting representations of $\mathfrak{o}(6)$ and $\mathfrak{su}(1,1)$ and
Howe duality takes place.

The maximal Abelian algebra of $\mathfrak{o}(6)$ is {$\os(2)\oplus\os(2)\oplus\os(2)$} and is generated by
the set \mbox{$\{L_{12},\,L_{34},\,L_{56}\}$}. The non-abelian part of its commutant in the representation of
$\mathcal{U}(\mathfrak{o}(6))$ on $\mathcal{H}$ is generated by the two invariants
\begin{align}
 K_1&=\frac{1}{8}\big({L_{12}}^2+{L_{34}}^2+{L_{13}}^2+{L_{23}}^2+{L_{14}}^2+{L_{24}}^2\big),\\
 K_2&=\frac{1}{8}\big({L_{34}}^2+{L_{56}}^2+{L_{35}}^2+{L_{36}}^2+{L_{45}}^2+{L_{46}}^2\big).\label{eq_K2}
\end{align}
Define $K_3$ by~ $[K_1,K_2]=K_3$. Working out the commutation relations of $K_3$ with $K_1$ and $K_2$, it is
found that they correspond to those \eqref{CR} of the Racah algebra
with the central parameters given by $d=-\frac{1}{8}\big(\mathcal{C}+{L_{12}}^2+{L_{34}}^2+{L_{56}}^2\big)$,\\
$e_1=-\frac{1}{64}\big(\mathcal{C}-{L_{12}}^2-4\big)\big({L_{34}}^2-{L_{56}}^2\big)$ and
$e_2=-\frac{1}{64}\big(\mathcal{C}-{L_{56}}^2-4\big)\big(L_{34}^2-L_{12}^2\big)$. For details see \cite{GVVZ}.
By abuse of notation we designate the abstract generators and their realizations by the same letter.

Regarding the $\mathfrak{su}(1,1)$ picture, let $J_\bullet^{(\mu,\nu, \rho, \lambda)} = J_\bullet^{(\mu)}
+J_\bullet^{(\nu)} + J_\bullet^{(\rho)} + J_\bullet^{ (\lambda)}$ denote the addition of the four metaplectic
representations labelled by the variables $\mu,\nu, \rho, \lambda$ all assumed different.  The corresponding
Casimir operator is $C^{(\mu,\nu, \rho, \lambda)}=(J_0^{(\mu,\nu, \rho, \lambda)})^2 $ $-J_+^{(\mu,\nu, \rho,
\lambda)}J_-^{(\mu,\nu, \rho, \lambda)}-J_0^{(\mu,\nu, \rho, \lambda)}$. Quite clearly, these actions of
$\mathfrak{su}(1,1)$ restricted to state vectors of four oscillators are paired with commuting actions of the
Lie algebra $\mathfrak{o}(4)$ of rotations in the four dimensions labelled by $\mu,\nu, \rho, \lambda$. It is
hence not surprising to find, owing to Howe duality, that  $C^{(1234)}=-2K_1$ and $C^{(3456)}=-2K_2$, namely
that the intermediate $\mathfrak{su}(1,1)$ Casimir operators corresponding to the recouplings of the first
four and last four of the six metaplectic representations are equal (up to a factor) to the Casimir elements
of the two corresponding $\mathfrak{o}(4)$ subalgebras of $\mathfrak{o}(6)$ which together generate as we
observed the non-trivial part of the commutant of \mbox{$\os(2)\oplus\os(2)\oplus\os(2)$} in
$\mathcal{U}(\mathfrak{o}(6))$. This entails the description of the Racah algebra as the commutant in
$\mathcal{U}(\mathfrak{su}(1,1)^{\otimes 3})$ of the action of $\mathfrak{su}(1,1)$ on $\mathcal{H}$.
Alternatively, picking the $\mathfrak{su}(1,1)$ representations associated to those of
\mbox{$\os(2)\oplus\os(2)\oplus\os(2)$} under Howe duality yields the sum of three irreducible representations
of  $\mathfrak{su}(1,1)$ belonging to the discrete series; these are realized as dynamical algebras of three
singular oscillators. Note that corresponding to the $\mathfrak{su}(1,1)$ representation
$J_\bullet^{(\mu,\nu)} = J_\bullet^{(\mu)} +J_\bullet^{(\nu)}$ is the Casimir $ C^{(\mu\nu)}
=-\frac{1}{4}\big({L_{\mu\nu}}^2+1\big)$. With the dependance on the polar angles ``rotated out", the total
Casimir element $C^{(123456)}$ becomes the Hamiltonian of the generic superintegrable system on the
two-sphere; the constants of motion are the quotiented intermediate Casimirs elements and the symmetry algebra
that they generate is hence that of Racah.

\section{More dual pictures -- an overview}\label{sec:more}

The main algebras of Askey--Wilson type have been studied recently from the commutant and Howe duality
viewpoints. We summarize in the following the main results and give in particular the dualities that are
involved.

\subsection{The Racah family}

The higher rank extension of the Racah algebra  defined as the algebra generated by all the intermediate
Casimir elements of $\sigma_1 \otimes \sigma_2 \otimes \dots \otimes \sigma_n (\Delta^{(n-1)}(\mathfrak{su}(1,
1)))$ can be described in the framework of the preceding section with the help of the dual pair
$(\mathfrak{o}(2n), \mathfrak{su}(1,1))$ using in this case the module formed by the state vectors of $2n$
harmonic oscillators. It is then seen to be dually the commutant of $\mathfrak{o}(2)^{\oplus n}$ in the
oscillator representation of $\mathcal{U}(\mathfrak{o}(2n))$ \cite{GVVZ2}.

The case $n=2$ is special and of particular interest since it pertains to the Clebsch--Gordan problem for
$\mathfrak{su}(1, 1)$, that is, the recoupling of the two irreducible representations $\sigma_1$ and
$\sigma_2$. There are no intermediate Casimirs here; the relevant operators associated to the direct product
basis and the recoupled one are respectively $M_1=\sigma_1(J_0) - \sigma_2(J_0)$ and the total Casimir ${M_2 =
(\sigma_1 \otimes \sigma_2) \Delta (C)}$. These are seen to obey the commutation relations of the Hahn algebra
\cite{FGVVZ}:
\begin{align}\label{eq_Hahn2}
\begin{aligned}{}
 [M_1,M_2] =M_3, \qquad
 [M_2,M_3]& =-2\{M_1,M_2\}+\delta_1,\\
 [M_3,M_1] &=-2{M_1}^{2}-4M_2+\delta_2,
\end{aligned}
\end{align}
where $\delta_1 = 4(\sigma_1(J_0) + \sigma_2(J_0))(\sigma_1(C) - \sigma_2(C))$  and $\delta_2 =
2(\sigma_1(J_0) + \sigma_2(J_0))^2 + (\sigma_1(C) + \sigma_2(C))$ are central. The name of the algebra comes
from the fact that the $3j$-coefficients involve dual Hahn polynomials. In the setup with four harmonic
oscillators, with $\mathcal{H}$ carrying the product of four metaplectic representations, Howe duality will
imply that the total Casimir element $C^{(1234)}$ of $\mathfrak{su}(1, 1)$ coincides with the Casimir of
$\mathfrak{o}(4)$ -- this is the same computation as the one described above. It is easily seen that
$\sigma_1(J_0) - \sigma_2(J_0)$ is derived from $\frac{1}{2}(N_1+N_2-N_3-N_4)$ under the quotient by
$\mathfrak{o}(2) \oplus \mathfrak{o}(2)$ with ${N_i = a_i^{\dagger} a_i}$, ~$i=1, \dots, 4$. It can in fact be
checked directly, again abusing notation, that ${M_1=\frac{1}{2}(N_1+N_2-N_3-N_4)}$ and ${M_2 =
-\frac{1}{4}\big({L_{12}}^2+{L_{34}}^2+{L_{13}}^2+{L_{23}}^2+{L_{14}}^2+{L_{24}}^2\big)}$ satisfy the
relations given in equation \eqref {eq_Hahn2} with $\delta_1=-\frac{1}{2}(N_1 + N_2 + N_3 + N_4 + 2)(
{L_{12}}^2 - {L_{34}}^2)$ and ${\delta_2 = \frac{1}{2}(N_1 + N_2 + N_3 + N_4 + 2)^2 - ( {L_{12}}^2 +
{L_{34}}^2 + 2)}$, in correspondance with the preceding expressions for $\delta_1$ and $\delta_2$ in the
realization $J_\bullet^{(1234)}$ of $\mathfrak{su}(1, 1)$. From the expressions of these last $M_1$ and $M_2$,
we can claim that the Hahn algebra is the commutant of $\mathfrak{o}(2) \oplus \mathfrak{o}(2)$ in
$\mathcal{U}(\mathfrak{u}(4))$ represented on $\mathcal{H}$.  Let us stress that it is the universal
enveloping algebra of $\mathfrak{u}(4)$ that intervenes here.

\subsection{The Bannai--Ito set}

The Bannai--Ito algebra \cite{DGTVZ} takes its name after the Bannai--Ito polynomials that enter in the Racah
coefficients of the Lie superalgebra $\mathfrak{osp}(1|2)$. This algebra has three generators $K_i,~
i=1,\dots, 3$ that satisfy the relations
\begin{align}\label{BIcr}
\begin{aligned}{}
 \{K_i, K_j\} = K_k + \omega_k, \qquad i \neq j \neq k \in \{1, 2, 3\}
\end{aligned}
\end{align}
with $\omega_i$ central and $\{X , Y\} = XY + YX$. The relevant reductive pair in this case is
$(\mathfrak{o}(6) , \mathfrak{osp}(1|2))$ and the representation space $\mathcal{H}$ is that of Dirac spinors
in six dimensions with the Clifford algebra generated by the elements $\gamma_{\mu}$ verifying
$\{\gamma_{\mu}, \gamma_{\nu}\} = -2\delta_{\mu \nu},~\mu, \nu = 1, \dots, 6$. That the pair
$(\mathfrak{o}(6) , \mathfrak{osp}(1|2))$ is dually represented on $\mathcal{H}$ is seen as follows: The
spinorial representation of $\mathfrak{o}(6)$ with generators
\begin{align}
\begin{aligned}{}
 J_{\mu \nu} = -i(x_{\mu}\partial_{\nu} - x_{\nu}\partial_{\mu})+ \Sigma_{\mu \nu}, \qquad
 \Sigma_{\mu \nu} = \frac{i}{2}\gamma_{\mu} \gamma_{\nu}
\end{aligned}
\end{align}
leaves invariant the following operators:
\begin{align}
\begin{aligned}{}
 J_- = -i\sum_{\mu=1}^6 \gamma_\mu \partial_\mu, \qquad J_+ = -i\sum_{\mu=1}^6 \gamma_\mu x_\mu, \qquad J_0 =
\sum_{\mu=1}^6 x_\mu \partial_\mu,
\end{aligned}
\end{align}
which in turn realize the commutation relations of the Lie superalgebra $\mathfrak{osp}(1|2)$: $[J_0, J_{\pm}]
= \pm J_{\pm}$,~ $\{ J_+, J_-\} = - 2J_0$ with $J_0$ even and $J_{\pm}$ odd. Howe duality thus takes place. As
a matter of fact, for any subset $A \subset \{1, \dots, 6\}$ of cardinality $|A|$ the operators
$J_-^{A}=-i\sum_{\mu\in A}\gamma_\mu\partial_\mu,~ J_+^{A}=-i\sum_{\mu\in A}\gamma_\mu x_\mu$, and
${J_0^{A}=\tfrac{|A|}{2}+\sum_{\mu\in A}x_\mu\partial_\mu}$ realize $\mathfrak{osp}(1|2)$. The Casimir element
of $\mathfrak{osp}(1|2)$ is given by $C = \frac{1}{2} ([J_- , J_+] - 1)S$ with $S$ the grade involution
obeying $S^2 = 1$, $[S , J_0] = 0$, $\{S , J_{\pm}\} = 0$. In the realizations at hand, $S^{A} = i^{|A|/2}
\prod_{\mu \in A} \gamma_{\mu}$ with $|A|$ even.

It can be checked that the operators
\begin{align*}
&\begin{aligned}
 K_1 = M_1 + \frac{3}{2} \Sigma_{12}\Sigma_{34}, \qquad K_2 = M_2 + \frac{3}{2} \Sigma_{34}\Sigma_{56},
 \qquad K_3 =  M_3 + \frac{3}{2} \Sigma_{12}\Sigma_{56},
\end{aligned}\\[.5em]
&\begin{aligned}
 \hspace{-1em}M_{1}&=\left(L_{12}\gamma_1\gamma_2+L_{13}\gamma_1\gamma_3+L_{14}\gamma_1\gamma_4
  +L_{23}\gamma_2\gamma_3+L_{24}\gamma_2\gamma_4+L_{34}\gamma_3\gamma_4\right)\Sigma_{12}\Sigma_{34},\\
 \hspace{-1em}M_{2}&=\left(L_{34}\gamma_3\gamma_4+L_{35}\gamma_3\gamma_5+L_{36}\gamma_3\gamma_6
  +L_{45}\gamma_4\gamma_5+L_{46}\gamma_4\gamma_6+L_{56}\gamma_5\gamma_6\right)\Sigma_{34}\Sigma_{56},\\
 \hspace{-1em}M_{3}&=\left(L_{12}\gamma_1\gamma_2+L_{15}\gamma_1\gamma_5+L_{16}\gamma_1\gamma_6
  +L_{25}\gamma_2\gamma_5+L_{26}\gamma_2\gamma_6+L_{56}\gamma_5\gamma_6\right)\Sigma_{12}\Sigma_{56}
\end{aligned}
\end{align*}
realize the relations \eqref{BIcr} of the Bannai-Ito algebra upon taking the following: ${\omega_{ij} =
2\Gamma_k\Gamma_{123} + 2\Gamma_i\Gamma_j}$, where ${\Gamma_1 = J_{12}}$, ${\Gamma_2 = J_{34}}$, ${\Gamma_3 =
J_{56}}$ and\linebreak ${\Gamma_{123} = \big(\frac{5}{2}-i\sum_{1\leq
\mu<\nu\leq6}L_{\mu\nu}\Sigma_{\mu\nu}\big)\Sigma_{12}\Sigma_{34}\Sigma_{56}}$. That these arise from dual
pictures is explained as follows (see \cite{GVVZ3} for details). On the one hand, $K_1, K_2, K_3$ are observed
to belong to the commutant in $\mathcal{U}(\mathfrak{o}(6))$ of the $\mathfrak{o}(2) \oplus \mathfrak{o}(2)
\oplus \mathfrak{o}(2)$ subalgebra of $\mathfrak{o}(6)$ spanned by $\{J_{12}, J_{34}, J_{56}\}$. On the other
hand, considering the Casimir elements $C^{A}$ of $\mathfrak{osp}(1|2)$ associated to the realization by the
operators $\{ J_0^{A}, J_{\pm}^{A}, S^{A}\}$, we find that $C^{(1234)} = K_1$,~ $C^{(3456)} = K_2$ and
$C^{(1256)} = K_3$. This confirms that the Bannai--Ito algebra can be dually presented either as the commutant
of  $\mathfrak{o}(2) \oplus \mathfrak{o}(2) \oplus \mathfrak{o}(2)$ in the spinorial representation of
$\mathcal{U}(\mathfrak{o}(6))$ or as the centralizer of the action of $\mathfrak{osp}(1|2)$ on $\mathcal{H}$.
These considerations can be extended to higher dimensions \cite{GVVZ3} so as to obtain analogously dual
commutant pictures for the Bannai--Ito algebras of higher ranks.

\subsection{The Askey--Wilson class}

The Askey--Wilson algebra can be presented as follows:
\begin{align}\label{eq:universalAW}
\begin{aligned}
 \frac{[K_A,K_B]_q}{q^{2}-q^{-2}}+K_C&=\frac{\gamma}{q+q^{-1}}, \qquad
 \frac{[K_B,K_C]_q}{q^{2}-q^{-2}}+K_A=\frac{\alpha}{q+q^{-1}},\\
 &\frac{[K_C,K_A]_q}{q^{2}-q^{-2}}+K_B=\frac{\beta}{q+q^{-1}},\\[.5em]
 \text{with}\qquad [A,B]_q&=qAB-q^{-1}BA \qquad \text{and\quad $\alpha$, $\beta$, $\gamma$ central}.
\end{aligned}
\end{align}
The $U_q(\su(1,1))$ algebra has three generators, $J_\pm$ and $J_0$, obeying ${[ J_0 \,, J_\pm ] = \pm J_\pm}$
and ${J_-J_+ - q^2 J_+J_- = q^{2J_0}\,[ 2J_0 ]_q}$ with ${[x]_q=\frac{q^{x}-q^{-x}}{q-q^{-1}}}$. Its coproduct
is defined by $\Delta (J_0) = J_0 \otimes 1 + 1 \otimes J_0$,~ $\Delta (J_\pm) = J_\pm \otimes q^{2J_0} + 1
\otimes J_\pm $. The Casimir operator $C$ of $U_q(\su(1,1))$ is given by ${C = J_+J_- q^{-2J_0+1} -
\frac{q}{(1-q^2)^2} \, \big( q^{2J_0-1} + q^{-2J_0+1} \big) + \frac{1+q^2}{(1-q^2)^2}}$.

The $q$-deformation $\mathfrak{o}_{q^{1/2}}(N)$ of $\mathfrak{o}(N)$ is defined as the algebra with generators
$L_{i,i+1}$ ($i=1,\dots,N-1$) obeying the relations
\begin{align*}
 & L_{i-1,i}\,L_{i,i+1}^2 - (q^{1/2}+q^{-1/2}) L_{i,i+1}\,L_{i-1,i}\,L_{i,i+1} + L_{i,i+1}^2\,L_{i-1,i} =
 -L_{i-1,i},\\
 & L_{i,i+1}\,L_{i-1,i}^2 - (q^{1/2}+q^{-1/2}) L_{i-1,i}\,L_{i,i+1}\,L_{i-1,i} + L_{i-1,i}^2\,L_{i,i+1} =
 -L_{i,i+1},\\
 & [ L_{i,i+1}, L_{j,j+1} ] = 0 \quad \text{for} \quad |i-j| > 1.
\end{align*}
We shall use the notation $L_{ik}^\pm = [L_{ij}^{\pm}\,,L_{jk}^{\pm}]_{q^{\pm 1/4}}$ for any $i<j<k$,
and by definition ${L_{i,i+1}^\pm=L_{i,i+1}}$.

The reductive pair $(\mathfrak{o}_{q^{1/2}}(6), U_q(\su(1,1))$ is the one which is of relevance for the
Askey--Wilson algebra. Let us indicate how $\mathfrak{o}_{q^{1/2}}(2n)$ and $U_q(\su(1,1))$ are dually
represented on the standard state space $\mathcal{H}$ of $2n$ independent $q$-oscillators described by
operators $\{A_i^\pm$, $A_i^0\}$ such that  ${[A_i^0, A_i^\pm] = \pm A_i^\pm}$,~ ${[A_i^-, A_i^+] =
q^{A_i^0}}$,\linebreak ${A_i^- A_i^+ - q A_i^+ A_i^- = 1}$,~ $i=1,\dots,2n$. The algebra $U_q(\su(1,1))$ is
represented on $\mathcal{H}$ by using the coproduct to embed it in the tensor product of $2n$ copies of the
$q$-deformation of the metaplectic representation, this gives
\begin{align}\label{qmet}
\begin{aligned}
 J_0^{(2n)}&=\Delta^{(2n-1)}\Big(\tfrac{1}{2}\left(A^{0}_i+\tfrac{1}{2}\right)\Big)
  =\frac{1}{2}\sum_{i=1}^{2n} \left(A_i^{0}+\frac{1}{2}\right),\\
 J_\pm^{(2n)}&=\Delta^{(2n-1)}\left(\frac{1}{[2]_{q^{1/2}}}(A^{\pm}_i)^{2}\right)
  =\frac{1}{[2]_{q^{1/2}}}\sum_{i=1}^{2n}
  \Bigg((A_i^{\pm})^{2}\prod_{j=i+1}^{2n}q^{A_j^{0}+\frac{1}{2}}\Bigg).
\end{aligned}
\end{align}
The algebra $\mathfrak{o}_{q^{1/2}}(2n)$ can also be realized in terms of $2n$ $q$-oscillators. The $2n-1$
generators take the form
\begin{align*}
 L_{i,i+1}= q^{-\frac{1}{2}(A_i^0+\frac{1}{2})} \Big( q^{\frac{1}{4}}A_i^+A_{i+1}^- -
 q^{-\frac{1}{4}}A_i^-A_{i+1}^+ \Big), \qquad i=1,\dots,2n-1.
\end{align*}
It can be checked that  $[J_0^{(2n)}, L_{i,i+1}]=[J_\pm^{(2n)}, L_{i,i+1}]=0,~ i=1,\dots,2n-1$, in other
words, that $U_q(\mathfrak{su}(1,1))$ and $\mathfrak{o}_{q^{1/2}}(2n)$ have commuting actions on the Hilbert
space $\mathcal{H}$ of $2n$ $q$-oscillators. This sets the stage for Howe duality. In order to connect with
the Askey--Wilson algebra we take $n=3$. The expressions of the operators $K_A$ and $K_B$ acting on
$\mathcal{H}$ that realize the relations \eqref{eq:universalAW} (together with the specific central elements)
are rather involved and we shall refer the reader to \cite{FGRV} for the formulas. We shall only stress that
these operators can be obtained in a dual way: They are affinely related to the generators of the commutant of
$\mathfrak{o}_{q^{1/2}}(2)^{\oplus 3}$ in $\mathfrak{o}_{q^{1/2}}(6)$ as well as to the intermediate
$U_q(\mathfrak{su}(1,1))$ Casimir elements $C^{(1234)} = \Delta^{(3)}(C) \otimes1 \otimes1$ and $C^{(3456)} =
1 \otimes1 \otimes \Delta^{(3)}(C)$ of the $q$-metaplectic representation (see \eqref{qmet}). This can be
extended to higher ranks by letting $n$ be arbitrary. For $n=2$ we are looking at the Clebsch--Gordan problen
for $U_q(\mathfrak{su}(1,1))$. The $q$-Hahn algebra that arises has two dual realizations \cite{FGRV2}: one as
the commutant of $\mathfrak{o}_{q^{1/2}}(2)^{\oplus 2}$ in $U_q(\mathfrak{u}(4))$ and the other in terms of
the following two $U_q(\mathfrak{su}(1,1))$ operators, $(\Delta(J_0)\otimes 1 \otimes 1) - (1\otimes 1 \otimes
\Delta(J_0)) $ and $\Delta^{(2)}(C)$ (the full Casimir element) in the $q$-metaplectic representation.

\section{Conclusion}\label{sec:conclusion}

This paper has offered a summary of how the quadratic algebras of Racah, Hahn, Bannai--Ito, Askey--Wilson and
$q$-Hahn types can be given dual descriptions as commutant of Lie algebras, superalgebras and quantum algebras.
The connection between these dual pictures is rooted in Howe dualities whose various expressions have been
stressed. The attentive reader will have noticed that the Clebsch--Gordan problem for $\mathfrak{osp}(1|2)$
has not been mentioned; this is because it has not been analyzed yet. We plan on adding this missing piece to
complete the picture.

\section*{Acknowledgments}

The authors thank Luc Frappat, Eric Ragoucy and Alexei Zhedanov for collaborations that led to the results
reviewed here.  JG holds an Alexander-Graham-Bell Scholarship from the Natural Science and Engineering
Research Council (NSERC) of Canada.  LV gratefully acknowledges his support from NSERC through a Discovery
Grant.

%
%

\end{document}